# Plankton-Oxygen Dynamics in the Context of Climate Change: A Fractional Model with A Probability Density Function Approach

Mahmoud M. El-Borai[1], Wagdy G. El-Sayed[2], Mahmoud A. Habib[3]

**Abstract**

Analyze how climate change affects marine oxygen production by modeling plankton–oxygen dynamics with a fractional-order nonlinear system and establishing rigorous conditions for the model's well-posedness.We formulate a three-dimensional system $d^\alpha x(t)/dt^\alpha = Ax(t)+f(x(t))$, where A is a diagonal matrix of order 3 and f is nonlinear. We (i) rigorously state the model, (ii) derive a Lipschitz constant for f under suitable assumptions, and (iii) prove existence, uniqueness, and continuous dependence on initial data using a fractional formula with a probability density kernel and a generalized Grönwall inequality.Under stated conditions, f satisfies a computable Lipschitz bound that yields existence and uniqueness of solutions for the fractional system. The solutions depend continuously on initial conditions, establishing well-posedness of the plankton–oxygen model.Introduces a fractional, PDF-kernel–based framework for plankton–oxygen dynamics and provides clean, general proofs of well-posedness via a generalized Grönwall approach, capturing memory effects that classical integer-order models miss.The results justify numerical simulation and sensitivity analyses of fractional marine-ecosystem models, providing a sound base for testing mitigation or management strategies affecting oxygen dynamics.Stronger theory for oxygen-cycle modeling can support evidence-based policies aimed at protecting marine ecosystems under global warming.

## Introduction

Over the past few decades, climate change has become a pressing issue [1]. It affects aquatic and terrestrial ecosystems through changes in temperature, moisture, atmospheric carbon dioxide, and solar UV radiation, affecting biodiversity, structure, function, and the ability to provide essential services [2]. It is a threat to almost all ecosystems and biological systems on Earth, and consequently to humans themselves, with 82% of the core ecological processes showing evidence of impact [3], causing a net negative effect on many ecosystem services, with losses exceeding gains in many places and affecting most people [4].

Oceans constitute approximately two-thirds of Earth's surface area and act as massive heat capacitors; thus, they are anticipated to be radically affected by climate change, as can be seen in [5], [6], [7]. That said, apart from the hydrophysical effects, oceans are massive ecosystems that are affected by global warming, which can be equally devastating as, or worse than, global flooding [8], [9]. It is known that plankton is the pillar of marine ecosystems and plays fundamental roles not only regrading food chains and global biogeochemical cycles [10], [11], but also when it comes to climate, composition of the atmosphere, and the amount of oxygen in it, as cited in [12]. In fact, plankton can be divided into two types: phytoplankton and zooplankton. Phytoplankton are plants that produce oxygen

[1] Department of Mathematics and Computer Science, Faculty of Science, Alexandria University, 21511 Alexandria, Egypt, Email: m_m_elborai@yahoo.com

[2] Department of Mathematics and Computer Science, Faculty of Science, Alexandria University, 21511 Alexandria, Egypt, Email: wagdygoma@alexu.edu.eg

[3] Department of Mathematics and Computer Science, Faculty of Science, Alexandria University, 21511 Alexandria, Egypt, Email: mahmoud.habib@alexu.edu.eg, (Corresponding Author)

via photosynthesis, like most plants, provided that sufficient light is available, whereas zooplankton are animals. The oxygen from phytoplankton is first released into the ocean and then into the atmosphere, contributing an estimated 70% of the atmospheric oxygen, as cited in [12]. As a result, we can predict that a decrease in the oxygen production of phytoplankton would have grave effects on life on Earth; thus, it is of great importance to study and identify potential threats to oxygen production [12]. It has been established that water temperature affects phytoplankton growth, respiration rate, and photosynthesis rate, and hence, net oxygen production is a function of temperature [12], which increases with global warming.

With what we established above, we study an important model provided in [12] using fractional calculus techniques. Fractional-order differential equations are superior tools for modeling complex phenomena because they capture their nonlocal and memory-dependent behaviors [13], [14], [15], [16], [17]. They are very powerful and provide a generalization of the integer-order derivative, which helps investigate the entire spectrum of values lying between two integer-order derivatives. Further applications can be found in [18], [19], [20], [21], [22], [23], [24], [25], [26], [27], [28], [29].

Although fractional-order derivatives have many definitions that can be found in [30], Caputo's definition provides the best framework to work with here. One vital reason is that it allows the use of classical initial conditions in a manner similar to that used in integer-order differential equations [31].

In this work, we analyze the mathematical model studying the dynamics of the marine ecosystem constituting of plankton and oxygen. The analysis is such that the model is a fractional-order differential equation, where the derivative is taken in Caputo's sense. In addition, in our analysis, we use the powerful formula given in [32], which provides elegance and simplicity. We organize the subsequent sections of this paper in the following manner: Section 2 presents some crucial groundwork material that will aid throughout the following sections. Section 3 provides the formulation of the fractional model and the necessary concepts and assumptions for it. In Section 4.1, we first prove the existence of a Lipschitz constant to the nonlinear part of the model $f$ and calculate it, which is significant for proving the existence of a solution to the model. We then prove this existence using the method of successive approximations. Next, in Section 4.2, we prove the uniqueness of the solution. Subsequently, in Section 4.3, we prove the continuous dependence on the initial condition which proves the wellposedness of the problem. The concluding remarks in Section 5 complete the article.

**Preliminaries**

Given a function $\phi$ with values in a Banach space, the integrals and derivatives that appear in the following definitions are taken in Bochner's sense.

**Definition 2.1.** The fractional integral of order $\beta > 0$ is given by

$$I^\beta \phi(t) = \frac{1}{\Gamma(\beta)} \int_0^t (t-s)^{\beta-1} \phi(s) \, ds.$$

**Definition 2.2.** If $0 < \beta \leq 1$, we define the derivative of order $\beta$ by

$$\frac{d^\beta \phi(t)}{dt^\beta} = \frac{1}{\Gamma(1-\beta)} \frac{d}{dt} \int_0^t \frac{\phi(s)}{(t-s)^\beta} \, ds.$$

**Definition 2.3.** (One-parameter Mittag-Leffler function)

Let $\gamma > 0$. Then the one-parameter Mittag-Leffler function of order $\gamma$ is defined by

$$E_\gamma(z) = \sum_{k=0}^\infty \frac{z^k}{\Gamma(k\gamma + 1)}.$$

We need the following generalized singular Grönwall inequality.

**Theorem 2.4.** [33] Let $J = [0, T], T < \infty$. Suppose $\beta > 0$, $h(t)$ is a nonnegative function locally integrable on $J$, and $q(t)$ is a bounded, nonnegative, non-decreasing continuous function defined on $J$. Let $p(t)$ be another nonnegative function, locally integrable on $J$, satisfying

$$p(t) \leq h(t) + q(t) \int_0^t (t-s)^{\beta-1} p(s) \, ds,$$

on this interval, then

$$p(t) \leq h(t) + \int_0^t \left[ \sum_{k=1}^{\infty} \frac{(q(t)\Gamma(\beta))^k}{\Gamma(k\beta)} (t-s)^{k\beta-1} h(s) \right] ds, \quad t \in J.$$

Setting $q(t) \equiv b$ in the theorem, we have the following inequality.

**Corollary 2.5.** [34] cf. [33]  Suppose $b \geq 0$, $\beta > 0$, and $h(t)$ is a nonnegative function locally integrable on $J$. Let $p(t)$ be nonnegative and locally integrable on $J$ with

$$p(t) \leq h(t) + b \int_0^t (t-s)^{\beta-1} p(s) \, ds,$$

on this interval, then

$$p(t) \leq h(t) + \int_0^t \left[ \sum_{k=1}^{\infty} \frac{(b\Gamma(\beta))^k}{\Gamma(k\beta)} (t-s)^{k\beta-1} h(s) \right] ds, \quad t \in J.$$

**Corollary 2.6.** [33]  Under the hypotheses of Theorem 2.4, if $h(t)$ is a non-decreasing function on $J$, then

$$p(t) \leq h(t) E_\beta\big(q(t)\Gamma(\beta) t^\beta\big).$$

**Formulation of the Problem**

Mathematical modeling using differential equations has proven to be very powerful for studying various natural phenomena. The marine ecosystem problem we deal with here is no exception. Based on the assumptions of competing species and the Lotka-Volterra model, as illustrated in Section 1, the interaction between oxygen concentration, phytoplankton concentration, and zooplankton concentration can be modeled. Phytoplankton produces oxygen via photosynthesis, which is necessary for its own respiration, as well as the zooplankton's. Zooplankton feeds on phytoplankton. The oxygen concentration may vary depending on biochemical reactions. The mathematical formula of the problem takes the form found in [12]

$$\begin{cases} \dfrac{dx_1}{dt} = \dfrac{Hc_0 x_2}{x_1 + c_0} - \dfrac{\delta x_1 x_2}{x_1 + c_2} - \dfrac{\nu x_3 x_1}{x_1 + c_3} - m x_1, \\ \dfrac{dx_2}{dt} = \left(\dfrac{Bx_1}{x_1 + c_1} - \gamma\right) x_2 - \dfrac{\beta x_2 x_3}{x_2 + h} - \sigma x_2, \\ \dfrac{dx_3}{dt} = \dfrac{\xi x_1^2}{x_1^2 + c_4^2} \cdot \dfrac{\beta x_2 x_3}{x_2 + h} - \mu x_3. \end{cases} \quad (1)$$

The meaning of each symbol is listed in Table 1. It should be noted that all symbols take non-negative values due to their biological meanings; otherwise, they would make no sense.

**Table 1. Description of Symbols.**

| Symbol | Description |
|---|---|
| $x_1(t)$ | Concentration of oxygen at time $t$ |
| $x_2(t)$ | Concentration of phytoplankton at time $t$ |
| $x_3(t)$ | Concentration of zooplankton at time $t$ |
| $c_0$ | Half-saturation constant for oxygen production |
| $c_1$ | Half-saturation constant for phytoplankton growth |
| $c_2$ | Half-saturation constant for phytoplankton respiration |
| $c_3$ | Half-saturation constant for zooplankton respiration |
| $c_4$ | Half-saturation constant related to zooplankton feeding efficiency |
| $h$ | Half-saturation constant for zooplankton predation |
| $H$ | Factor accounting for the environmental effect on oxygen production in phytoplankton |
| $\delta$ | Maximum per capita phytoplankton respiration rate |
| $\nu$ | Maximum per capita zooplankton respiration rate |
| $B$ | Maximum phytoplankton per capita growth rate in presence of high oxygen concentration |
| $\beta$ | Maximum predation rate of zooplankton on phytoplankton |

| $\xi$ | Maximum feeding efficiency of zooplankton |
| --- | --- |
| $m$ | Rate of oxygen loss due to natural depletion (e.g., biochemical reactions) |
| $\gamma$ | Intensity of intraspecific competition among phytoplankton |
| $\sigma$ | Natural mortality rate of phytoplankton |
| $\mu$ | Natural mortality rate of zooplankton |

We now formulate the fractional version of System (1).

Consider the Banach space $(\mathbb{R}^3, \|\cdot\|_1)$. where $\|\cdot\|_1$ is the usual 1-norm. Take $J = [0, T]$, $T < \infty$. Define the nonnegative real line as $\mathbb{R}_{\geq 0} = \{y \in \mathbb{R} \mid y \geq 0\}$, and the nonnegative octant in $\mathbb{R}^3$ as $\mathbb{R}^3_{\geq 0} = \{z = [z_1, z_2, z_3]^T \in \mathbb{R}^3 \mid z_i \geq 0, i = 1, 2, 3\}$. Now, for $i = 1, 2, 3$, we may define the following functions

$x_i: J \to \Omega_i \subset R_{\geq 0}$,

$x: J \to \Omega \subset R^3_{\geq 0}$, where $\Omega = \prod_{i=1}^3 \Omega_i$,

$f: \Omega \to R^3_{\geq 0}$.

Additionally, we define $\|x(t)\| = max_{t \in J} \|x(t)\|_1 = max_{t \in J} \sum_{i=1}^3 |x_i(t)|$ and $\|f(x(t))\| = max_{t \in J} \|f(x(t))\|_1 = max_{t \in J} \sum_{i=1}^3 |f_i(x(t))|$. Also, for the closed linear matrix $A$, we have $\|A\|_1 = \sum_{i,j=1}^3 |a_{i,j}| = m + \sigma + \mu$.

Define the fractional-order evolution equation with order $0 < \alpha \leq 1$ in the Banach space $(\mathbb{R}^3, \|\cdot\|_1)$ as

$$\begin{cases} \dfrac{d^\alpha x(t)}{dt^\alpha} = Ax(t) + f(x(t)), \\ x(0) = x_0, \end{cases} \qquad (2)$$

where $x(t) = [x_1(t), x_2(t), x_3(t)]^T$, $A = diag(-m, -\sigma, -\mu)$, and

$$f(x(t)) = \begin{bmatrix} \dfrac{Hc_0 x_2(t)}{x_1(t) + c_0} - \dfrac{\delta x_2(t) x_1(t)}{x_1(t) + c_2} - \dfrac{\nu x_1(t) x_3(t)}{x_1(t) + c_3} \\ \left(\dfrac{B x_1(t)}{x_1(t) + c_1} - \gamma x_2(t)\right) x_2(t) - \dfrac{\beta x_2(t) x_3(t)}{x_2(t) + h} \\ \dfrac{\xi \beta x_1^2(t) x_2(t) x_3(t)}{(x_1^2(t) + c_4^2)(x_2(t) + h)} \end{bmatrix}.$$

Following the work in [32], we see that Equation (2) is equivalent to

$$\begin{cases} x(t) = \displaystyle\int_0^\infty \zeta_\alpha(\theta) e^{t^\alpha \theta A} x_0 \, d\theta + \alpha \int_0^t \int_0^\infty \theta(t-\eta)^{\alpha-1} \zeta_\alpha(\theta) e^{(t-\eta)^\alpha \theta A} f(x(\eta)) \, d\theta \, d\eta, \\ x(0) = x_0, \end{cases} \qquad (3)$$

which is the formula we use in our analysis. Here, $\zeta_\alpha(\theta)$ is a probability density function defined on $(0, \infty)$.

**Analysis of the Problem.**

**Existence of Solutions.**

**Lemma 4.1.1.** The function $f$ in the System (2) is Lipschitz with constant $L$.

*Proof.* First, since each component $x_i(t)$ of $x(t)$ is nonnegative and bounded on $J$, then we have the following

(C1) there exist constants $M_i$, such that $0 \leq x_i(t) \leq M_i$, on $J$, for $i = 1, 2, 3$,

(C2) $x_1 + c_k \geq c_k$, for $k = 0, 1, 2, 3$. Similarly, $x_2 + h \geq h$ and $x_1^2 + c_4^2 \geq c_4^2$.

Let $f = [f_1, f_2, f_3]^T$, $f_1 = F_1 + F_2 + F_3$. For $F_1(x) = \dfrac{Hc_0 x_2}{x_1 + c_0}$, we have

$$|F_1(x) - F_1(y)| = Hc_0 \left| \dfrac{x_2}{x_1 + c_0} - \dfrac{y_2}{y_1 + c_0} \right|$$

$$= Hc_0 \left| \dfrac{x_2(y_1 + c_0) - y_2(x_1 + c_0)}{(x_1 + c_0)(y_1 + c_0)} \right|,$$

1049

thus, and using (C1) and (C2), we get

$$|F_1(x) - F_1(y)| \leq \frac{H}{c_0}[M_2|x_1 - y_1| + (M_1 + c_0)|x_2 - y_2|].$$

Set

$$K_1 = \frac{HM_2}{c_0}, \quad K_2 = \frac{H(M_1 + c_0)}{c_0}.$$

For the second term $F_2(x) = \frac{-\delta x_1 x_2}{x_1 + c_2}$, we have

$$|F_2(x) - F_2(y)| = \delta \left| \frac{x_1 x_2}{x_1 + c_2} - \frac{y_1 y_2}{y_1 + c_2} \right|$$
$$= \delta \left| \frac{x_1 x_2 (y_1 + c_2) - y_1 y_2 (x_1 + c_2)}{(x_1 + c_2)(y_1 + c_2)} \right|,$$

thus, in a similar manner, we get

$$|F_2(x) - F_2(y)| \leq \frac{\delta}{c_2^2}[M_2 c_2 |x_1 - y_1| + M_1(M_1 + c_2)|x_2 - y_2|].$$

Set

$$K_3 = \frac{\delta M_2}{c_2}, \quad K_4 = \frac{\delta M_1(M_1 + c_2)}{c_2^2}.$$

For the third term $F_3(x) = \frac{-\nu x_1 x_3}{x_1 + c_3}$, we have, by a similar argument,

$$|F_3(x) - F_3(y)| \leq \frac{\nu}{c_3^2}[M_3 c_3 |x_1 - y_1| + M_1(M_1 + c_3)|x_3 - y_3|].$$

Set

$$K_5 = \frac{\nu M_3}{c_3}, \quad K_6 = \frac{\nu M_1(M_1 + c_3)}{c_3^2}.$$

Thus

$$|f_1(x) - f_1(y)| \leq (K_1 + K_3 + K_5)|x_1 - y_1| + (K_2 + K_4)|x_2 - y_2| + K_6|x_3 - y_3|.$$

Now, let $f_2 = G_1 + G_2$. For $G_1(x) = \left(\frac{Bx_1}{x_1 + c_1} - \gamma x_2\right) x_2$, we have

$$|G_1(x) - G_1(y)| = \left| \left(\frac{Bx_1}{x_1 + c_1} - \gamma x_2\right) x_2 - \left(\frac{By_1}{y_1 + c_1} - \gamma y_2\right) y_2 \right|$$
$$= \left| B\left(\frac{x_1 x_2}{x_1 + c_1} - \frac{y_1 y_2}{y_1 + c_1}\right) - \gamma(x_2^2 - y_2^2) \right|.$$

Simplifying with the aid of (C1) and (C2), we get

$$|G_1(x) - G_1(y)| \leq \left| B\left[\frac{x_1}{x_1 + c_1}(x_2 - y_2) + \frac{c_1 y_2}{(x_1 + c_1)(y_1 + c_1)}(x_1 - y_1)\right] - \gamma(x_2 + y_2)(x_2 - y_2) \right|$$
$$\leq B\left[\frac{M_1}{c_1}|x_2 - y_2| + \frac{c_1 M_2}{c_1^2}|x_1 - y_1|\right] + 2M_2 \gamma |x_2 - y_2|.$$

Set

$$K_7 = \frac{BM_2}{c_1}, \quad K_8 = \frac{BM_1}{c_1} + 2M_2.$$

For $G_2(x) = -\beta \frac{x_2 x_3}{x_2 + h}$, we have

$$|G_2(x) - G_2(y)| = \beta \left| \frac{x_2 x_3}{x_2 + h} - \frac{y_2 y_3}{y_2 + h} \right|$$

Adding and subtracting $\frac{x_2 y_3}{x_2+h}$, using (C1) and (C2), and collecting terms in a similar way as precedes, we get

$$|G_2(x) - G_2(y)| \leq \beta \left| \frac{h y_3}{(x_2+h)(y_2+h)}(x_2 - y_2) + \frac{x_2}{x_2+h}(x_3 - y_3) \right|$$
$$\leq \beta \left[ \frac{h M_3}{h^2}|x_2 - y_2| + \frac{M_2}{h}|x_3 - y_3| \right].$$

Set

$$K_9 = \frac{\beta M_3}{h}, \qquad K_{10} = \frac{\beta M_2}{h}.$$

Thus

$$|f_2(x) - f_2(y)| \leq K_7|x_1 - y_1| + (K_8 + K_9)|x_2 - y_2| + K_{10}|x_3 - y_3|.$$

Let $f_3(x) = \frac{\xi \beta x_1^2 x_2 x_3}{(x_1^2 + c_4^2)(x_2 + h)}$. We have

$$|f_3(x) - f_3(y)| = \xi \beta \left| \frac{x_1^2 x_2 x_3}{(x_1^2 + c_4^2)(x_2 + h)} - \frac{y_1^2 y_2 y_3}{(y_1^2 + c_4^2)(y_2 + h)} \right|,$$

adding and subtracting both $\frac{y_1^2 x_2 x_3}{(y_1^2 + c_4^2)(x_2 + h)}$ and $\frac{y_1^2 y_2 x_3}{(y_1^2 + c_4^2)(y_2 + h)}$, then simplifying, we get

$$|f_3(x) - f_3(y)| = \xi \beta \left| \frac{c_4^2 x_2 x_3 (x_1 + y_1)}{(x_2 + h)(x_1^2 + c_4^2)(y_1^2 + c_4^2)}(x_1 - y_1) + \frac{h y_1^2}{(y_1^2 + c_4^2)(x_2 + h)(y_2 + h)}(x_2 - y_2) \right.$$
$$\left. + \frac{y_1^2 y_2}{(y_1^2 + c_4^2)(y_2 + h)}(x_3 - y_3) \right|.$$

Thus, via (C1) and (C2), we obtain

$$|f_3(x) - f_3(y)| \leq \xi \beta \left[ \frac{2 c_4^2 M_1 M_2 M_3}{h c_4^4}|x_1 - y_1| + \frac{h M_1^2}{h^2 c_4^2}|x_2 - y_2| + \frac{M_1^2 M_2}{h c_4^2}|x_3 - y_3| \right].$$

Set

$$K_{11} = \frac{2 \xi \beta M_1 M_2 M_3}{h c_4^2}, \qquad K_{12} = \frac{\xi \beta M_1^2}{h c_4^2}, \qquad K_{13} = \frac{\xi \beta M_1^2 M_2}{h c_4^2}.$$

Thus

$$|f_3(x) - f_3(y)| \leq K_{11}|x_1 - y_1| + K_{12}|x_2 - y_2| + K_{13}|x_3 - y_3|.$$

Now, we have

$$\|f(x) - f(y)\|_1 = |f(x_1) - f(y_1)| + |f(x_2) - f(y_2)| + |f(x_3) - f(y_3)|$$
$$\leq (K_1 + K_3 + K_5 + K_7 + K_{11})|x_1 - y_1| + (K_2 + K_4 + K_8 + K_9 + K_{12})|x_2 - y_2|$$
$$+ (K_6 + K_{10} + K_{13})|x_3 - y_3|$$
$$= L_1|x_1 - y_1| + L_2|x_2 - y_2| + L_3|x_3 - y_3|$$
$$\leq \max\{L_1, L_2, L_3\} \|x - y\|_1.$$

Taking $\max_{t \in J}$ for both sides and setting $L = \max\{L_1, L_2, L_3\}$, we have $\|f(x) - f(y)\| \leq L \|x - y\|$. ∎

**Theorem 4.1.2.** Equation (3) has a solution.

*Proof.* Using Picard's method of successive approximations in Equation (3), set

$$x_{n+1}(t) = \int_0^\infty \zeta_\alpha(\theta) e^{t^\alpha \theta A} x_0 \, d\theta + \alpha \int_0^t \int_0^\infty \theta(t - \eta)^{\alpha - 1} \zeta_\alpha(\theta) e^{(t-\eta)^\alpha \theta A} f(x_n(\eta)) \, d\theta \, d\eta,$$

then

$$x_{n+1}(t) - x_n(t) = \alpha \int_0^t \int_0^\infty \theta(t - \eta)^{\alpha - 1} \zeta_\alpha(\theta) e^{(t-\eta)^\alpha \theta A} \cdot [f(x_n(\eta)) - f(x_{n-1}(\eta))] \, d\theta \, d\eta.$$

Taking the norm on both sides and using Lemma 4.1, we have

$$\|x_{n+1}(t) - x_n(t)\| \leq \alpha \int_0^t \int_0^\infty \theta(t-\eta)^{\alpha-1} \zeta_\alpha(\theta) \|e^{(t-\eta)^\alpha \theta A}\| \|f(x_n(\eta)) - f(x_{n-1}(\eta))\| \, d\theta \, d\eta$$

$$\leq \alpha L \int_0^t \int_0^\infty \theta(t-\eta)^{\alpha-1} \zeta_\alpha(\theta) \|e^{(t-\eta)^\alpha \theta A}\| \|x_n(\eta) - x_{n-1}(\eta)\| \, d\theta \, d\eta.$$

Then, following El-Borai in [32], we have

$$\|x_2(t) - x_1(t)\| \leq \frac{Mt^\alpha}{\alpha},$$

where $M = \alpha K L \sup_{t,\theta} \|e^{t^\alpha \theta A} x_0\|$, $t \in J$, $\theta \in [0, \infty)$, $L$ is the Lipschitz constant of $f(x(t))$, and $K$ is a constant.

Hence, by induction, we get

$$\|x_{n+1}(t) - x_n(t)\| \leq \frac{M^n t^{n\alpha} (\Gamma(\alpha))^n}{\Gamma(n\alpha + 1)}.$$

Thus, the series $\sum_{k=1}^\infty [x_k(t) - x_{k-1}(t)]$ converges uniformly on $J$ with respect to the norm to a continuous function $x$.

Letting

$$x(t) = x_0 + \sum_{k=1}^\infty [x_k(t) - x_{k-1}(t)], \quad t \in J,$$

implies that a solution $x(t)$ exists. ∎

**Uniqueness of the solution.**

**Theorem 4.2.** The solution to Equation (3) is unique.

*Proof.* Let $x(t)$ and $y(t)$ be two solutions to Equation (3), then, and by using Lemma 4.1, we get

$$\|x(t) - y(t)\| \leq \alpha \int_0^t \int_0^\infty \theta(t-\eta)^{\alpha-1} \zeta_\alpha(\theta) \|e^{(t-\eta)^\alpha \theta A}\| \|f(x(\eta)) - f(y(\eta))\| \, d\theta \, d\eta$$

$$\leq \alpha L \int_0^t \int_0^\infty \theta(t-\eta)^{\alpha-1} \zeta_\alpha(\theta) \|e^{(t-\eta)^\alpha \theta A}\| \|x(\eta) - y(\eta)\| \, d\theta \, d\eta$$

$$\leq M \int_0^t (t-\eta)^{\alpha-1} \|x(\eta) - y(\eta)\| \, d\eta,$$

where $M$ is the same as in the proof of Theorem 4.1.2.

Setting $\|x(t) - y(t)\| = U(t)$, we have

$$U(t) \leq M \int_0^t (t-\eta)^{\alpha-1} U(\eta) \, d\eta.$$

Using Corollary 2.5, we obtain $U(t) \leq 0$ thus, $x(t) = y(t)$. ∎

**Wellposedness: Continuous Dependence on the Initial Condition.**

**Theorem 4.3.** The solution continuously depends on the initial conditions.

*Proof.* Using Equation (2), we have

$$\begin{cases} \dfrac{d^\alpha x_{n(t)}}{dt^\alpha} = A x_n(t) + f(x_n(t)), \\ \quad x(0) = a_n. \end{cases}$$

We show that $a_n \to a$ implies $x_n(t) \to x(t)$ as $n \to \infty$.

The convergence of the sequence $(a_n)$ implies that

(D1) for every $\varepsilon > 0$, there exists a natural $N > 0$, such that $n \geq N$ implies $\|a_n - a\| \leq \varepsilon$.

Using Equation (3), we have

$$x_n(t) = \int_0^\infty \zeta_\alpha(\theta) e^{t^\alpha \theta A} a_n \, d\theta + \alpha \int_0^t \int_0^\infty \theta(t-\eta)^{\alpha-1} \zeta_\alpha(\theta) e^{(t-\eta)^\alpha \theta A} f(x_n(\eta)) \, d\theta d\eta. \tag{4}$$

Subtracting Equation (3) from Equation (4), taking the norm, and using Lemma 4.1, we obtain

$$\|x_n(t) - x(t)\| \leq \int_0^\infty \zeta_\alpha(\theta) \|e^{t^\alpha \theta A}\| \|a_n - a\| \, d\theta$$

$$+ \alpha \int_0^t \int_0^\infty \theta(t-\eta)^{\alpha-1} \zeta_\alpha(\theta) \|e^{(t-\eta)^\alpha \theta A}\| \|f(x_n(\eta)) - f(x(\eta))\| \, d\theta \, d\eta$$

$$\leq K^* \|a_n - a\| + \alpha L \int_0^t \int_0^\infty \theta(t-\eta)^{\alpha-1} \zeta_\alpha(\theta) \|e^{(t-\eta)^\alpha \theta A}\| \|x_n(\eta) - x(\eta)\| \, d\theta \, d\eta.$$

Using (D1), we obtain

$$\|x_n(t) - x(t)\| \leq \varepsilon K^* + M \int_0^t (t-\eta)^{\alpha-1} \|x_n(\eta) - x(\eta)\| \, d\eta,$$

where $M$ is the same as in the proof of Theorem 4.1.2.

Setting $\|x_n(t) - x(t)\| = V_n(t)$, then

$$V_n(t) \leq \varepsilon K^* + M \int_0^t (t-\eta)^{\alpha-1} V_n(\eta) \, d\eta.$$

Using Corollary 2.5 and Corollary 2.6, we get

$$V_n(t) \leq \varepsilon K^* E_\alpha(M \Gamma(\alpha) t^\alpha),$$

for every $\varepsilon$; hence the required. ∎

## Conclusion

Since the most part of the oxygen in the atmosphere is produced by plankton in oceans, and to obtain a deep understanding of the climate change effect on this system, this paper studied a mathematical model governing the plankton-oxygen dynamics using a fractional-order differential equation as a generalization of its first-order counterpart. We proved the existence and uniqueness of the solution, as well as the wellposedness of problem via continuous dependence on the initial conditions, using a fractional formula with a probability density function, and a modified Grönwall inequality. In contrast to previous models, our approach integrates memory effects through the fractional derivative, making it more suitable for long-term environmental forecasting. Moreover, the elegant and powerful formula we employed simplifies the proofs. This work can be regarded as a theoretical foundation for more extensive simulations in climate modeling.

## Declarations

### Conflict of Interest

The authors declare that there is no conflict of interest.

### Ethics Approval Statement

This study did not involve any human participants or animals and thus did not require ethical approval.

### Credit Authorship Contribution Statement

*Mahmoud M. El-Borai:* conceptualization, formal analysis, investigation, methodology, resources, supervision

*Wagdy G. El-Sayed:* formal analysis, investigation, supervision, writing - review & editing

*Mahmoud A. Habib:* formal analysis, resources, writing - original draft, writing - review & editing

### Data Access Statement

Not applicable.

**Funding Sources**

This research did not receive any specific grant from funding agencies in the public, commercial, or not-for-profit sectors.

## Acknowledgements

We express our sincere gratitude to our colleagues who generously provided feedback and encouragement throughout the development of this work.